\documentclass[a4paper,12pt,leqno]{amsart}
\usepackage{amsmath,amssymb,mathrsfs}
\usepackage[matrix,arrow]{xy} 
\usepackage{xcolor,hyperref}
\hypersetup{colorlinks=true,citecolor=black,linkcolor=black}
\newcommand{\h}{\hbox}
\newcommand{\q}{\quad}
\newcommand{\nin}{\noindent}

\newcommand{\ms}{\par\medskip}
\newcommand{\sk}{\par\smallskip}

\newcommand{\msn}{\par\medskip\noindent}

\newcommand{\ges}{\geqslant}
\newcommand{\les}{\leqslant}
\newcommand{\one}{\hskip1pt}

\newcommand{\mcup}{\hbox{$\bigcup$}}
\newcommand{\msum}{\hbox{$\sum$}}
\newcommand{\mopl}{\hbox{$\bigoplus$}}
\newcommand{\mprod}{\hbox{$\prod$}}
\newcommand{\D}{{\mathcal D}}
\newcommand{\Hc}{{\mathcal H}}
\newcommand{\I}{{\mathcal I}}
\newcommand{\Lc}{{\mathcal L}}
\newcommand{\M}{{\mathcal M}}
\newcommand{\OO}{{\mathcal O}}

\newcommand{\PP}{{\mathbb P}}
\newcommand{\Q}{{\mathbb Q}}
\newcommand{\C}{{\mathbb C}}
\newcommand{\N}{{\mathbb N}}
\newcommand{\R}{{\mathbb R}}

\newcommand{\Z}{{\mathbb Z}}
\newcommand{\ft}{\widetilde{f}}
\newcommand{\iti}{\widetilde{i}}
\newcommand{\jt}{\widetilde{j}}
\newcommand{\Ct}{\widetilde{C}}
\newcommand{\Dt}{\widetilde{D}}
\newcommand{\Et}{\widetilde{E}}
\newcommand{\Ht}{\widetilde{H}}
\newcommand{\Xt}{\widetilde{X}}
\newcommand{\Yt}{\widetilde{Y}}

\newcommand{\pit}{\widetilde{\pi}}
\newcommand{\rhot}{\widetilde{\rho}}
\newcommand{\Lct}{\widetilde{\mathcal L}}
\newcommand{\al}{\alpha}
\newcommand{\be}{\beta}
\newcommand{\ga}{\gamma}
\newcommand{\Ga}{\Gamma}
\newcommand{\de}{\delta}

\newcommand{\la}{\lambda}

\newcommand{\si}{\sigma}

\newcommand{\ep}{\varepsilon}

\newcommand{\Om}{\Omega}
\newcommand{\Cr}{C_{\rm red}}
\newcommand{\Ctr}{\widetilde{C}_{\rm red}}
\newcommand{\YS}{Y_{\mhs S}}
\newcommand{\YtS}{\widetilde{Y}_{\mhs S}}
\newcommand{\ddd}{{\rm d}}
\newcommand{\Gr}{{\rm Gr}}
\newcommand{\Sp}{{\rm Sp}}
\newcommand{\Spf}{{\rm Sp}_{\mhs f}}

\newcommand{\tos}{\,{\to}\,}
\newcommand{\eq}{\,{=}\,}
\newcommand{\defs}{\,{:=}\,}
\newcommand{\nes}{\,{\ne}\,}
\newcommand{\ins}{\,{\in}\,}
\newcommand{\caps}{\,{\cap}\,}
\newcommand{\sst}{\,{\subset}\,}
\newcommand{\stm}{\,{\setminus}\,}
\newcommand{\gess}{\,{\ges}\,}
\newcommand{\less}{\,{\les}\,}
\newcommand{\sgt}{\,{>}\,}
\newcommand{\slt}{\,{<}\,}
\newcommand{\col}{\,{:}\,}
\newcommand{\pl}{\one {+}\one}
\newcommand{\mi}{\one {-}\one}
\newcommand{\bl}{\bigl}
\newcommand{\br}{\bigr}

\newcommand{\ssb}{\raise.15ex\h{${\scriptscriptstyle\bullet}$}}
\newcommand{\ssc}{\,\raise.15ex\h{${\scriptstyle\circ}$}\,}
\newcommand{\onto}{\twoheadrightarrow}
\newcommand{\into}{\hookrightarrow}
\newcommand{\simto}{\,\,\rlap{\hskip1.5mm\raise1.4mm\hbox{$\sim$}}\hbox{$\longrightarrow$}\,\,}
\newcommand{\mhs}{\hskip-1pt}

\makeatletter
\renewcommand\section{\@startsection{section}{1}{0pt}{-2ex plus -1ex minus -.2ex}{-2.3ex plus-.2ex}{\centering\normalfont\bfseries}}
\makeatother
\theoremstyle{plain}
\newtheorem{thm}{Theorem}[section]

\newtheorem{ithm}{Theorem}
\newtheorem{icor}{Corollary}

\theoremstyle{definition}
\newtheorem{rem}[thm]{Remark}

\newtheorem{irem}{Remark}

\begin{document}
\title[Spectrum of cones of hypersurfaces]{Spectrum of cones of projective hypersurfaces\\with singularities isolated}
\author[S.-J. Jung]{Seung-Jo Jung}
\address{S.-J. Jung : Department of Mathematics Education, and Institute of Pure and Applied Mathematics, Jeonbuk National University, Jeonju, 54896, Korea}
\email{seungjo@jbnu.ac.kr}
\author[M. Saito]{Morihiko Saito}
\address{M. Saito : RIMS Kyoto University, Kyoto 606-8502 Japan}
\email{msaito@kurims.kyoto-u.ac.jp}
\author[Y. Yoon]{Youngho Yoon}
\address{Y. Yoon : Department of Mathematics, Chungbuk National University, Cheongju-si, Chungcheongbuk-do, 28644, Korea}
\email{mathyyoon@gmail.com}
\thanks{This work was partially supported by the research grant of the Chungbuk National University in 2022 and by National Research Foundation of Korea (the first author: NRF-2021R1C1C1004097 and the third author: RS-2023-00245670).}
\begin{abstract} Let $Z$ be a projective hypersurface such that its underlying reduced variety has only isolated singularities. In case its irreducible components have constant multiplicities, for instance if $\dim Z>1$, we show that the spectrum of its cone can be described by using the spectral numbers at singular points of the reduced hypersurface and the global degree. In the non-reduced plane curve case, assuming that the underlying reduced curve has only semi-weighted-homogeneous singularities, we express the spectrum of the cone in terms of  the local weights and the weighted degrees and multiplicities of local irreducible components together with the degrees and multiplicities of global ones. These generalize a formula for reduced line arrangements. In the non-reduced ordinary (that is, semi-homogeneous) singularity case, the second formula is essentially equivalent to the one obtained by the third-named author.
\end{abstract}
\maketitle
\part*{Introduction} \label{intro}
\nin
The spectrum and Bernstein-Sato polynomial are quite interesting invariants of hypersurface singularities, although they are not necessarily easy to calculate in general even in the case of the cone of a projective hypersurface $Z\sst Y\defs\PP^n$. Let $f$ be a defining homogeneous polynomial of $Z$. Let $\Spf(t)\eq\msum_{\al\in\Q}\,n_{f,\al}t^{\al}$ be the spectrum of $f$ at 0, see Section\,\,\ref{S1} below. Put $d\defs\deg f\eq\deg Z$. Note that $n_{f,\al}\eq0$ if $\al d\,{\notin}\,\Z\cap(0,(n{+}1)d)$, see \eqref{1.3} below.
\sk
We assume in this paper that the reduced variety $Z'$ underlying $Z$ has only {\it isolated\one} singularities. Let $p_j$ ($j\ins[1,q'']$) be the singular points of $Z'$ with $\al_{j,l}$ ($l\ins[1,\mu_j]$) the spectral numbers of $(Z',p_j)$ counted with multiplicities, where $\mu_j$ is the Milnor number at $p_j$, see \cite{St0} (and Section\,\,\ref{S1} below). Set
$$M_j(\be):=\#\bl\{l\ins[1,\mu_j]\,\big|\,\be{-}1\less\al_{j,l}\slt\be\br\}.$$
Let $f'$ be a defining polynomial of $Z'$ with $d'\defs\deg f'$. Let $Z_k$ ($k\ins[1,r]$) be the irreducible components of $Z'$. Consider the condition that there is an integer $m$ such that
\msn
(CM)\hfill${\rm mult}_{Z_k}Z\eq m\,\,\,\,(\forall\,k\ins[1,r]),\,$ that is, $\,f\eq f'{}^m$.\q{}\hfill{}
\msn
This is always satisfied if $n\sgt2$. Define the integers $\ga_i$ by $\msum_{i\in\Z}\,\ga_i\one t^i\eq(t+\cdots+t^{d'-1})^{n+1}$. We can show the following.

\begin{ithm} \label{T1}
For any integers $i\ins[1,(n{+}1)d'{-}1]$, there are equalities
\begin{equation} \label{1}
n_{f',\frac{i}{d'}}=\ga_i-\msum_{j=1}^{q''}\,M_j\bl(\tfrac{i}{d'}\br).
\end{equation}
If $n\sgt2$ or more generally, condition {\rm(CM)} is satisfied, we have the following equalities for any integers $i\ins[1,d']$, $l\ins[0,m{-}1]$, $p'\in[0,n]:$
\begin{equation} \label{2}
n_{f,\frac{i}{d}+\frac{l}{m}+p'}=\begin{cases}n_{f',\frac{i}{d'}+p'}+(-1)^n&\h{if}\,\,\,i\eq d',\,p'\eq n,\, l\ne m{-}1,\\n_{f',\frac{i}{d'}+p'}&\h{otherwise.}\end{cases}
\end{equation}
\end{ithm}

Using Remark\,\,\ref{R4} and \eqref{6.1}, \eqref{6.2} below (and $\tbinom{-1}{2}\eq1$ for the case $\tfrac{i}{d}{+}2\eq3$), this implies the following.

\begin{icor} \label{C1}
If $Z$ is reduced and $n\eq2$, we have for integers $i\ins[1,d]$ the equalities
\begin{equation} \label{3}
\aligned n_{f,\frac{i}{d}}&=\tbinom{i-1}{2}-\msum_{j=1}^{q''}\,M_j\bl(\tfrac{i}{d}\br),\\ n_{f,\frac{i}{d}+1}&=(i{-}1)(d{-}i{-}1)\pl\tbinom{d}{2}-\msum_{j=1}^{q''}M_j\bl(\tfrac{i}{d}\pl1\br),\\ n_{f,\frac{i}{d}+2}&=\tbinom{d-i-1}{2}-\msum_{j=1}^{q''}\,M_j\bl(\tfrac{i}{d}\pl2\br)-\de_{i,d}.\endaligned
\end{equation}
\end{icor}

This is essentially a generalization of a formula for reduced line arrangements \cite{BS} via the argument in Section\,\,\ref{S6} below.
\sk
From now on we {\it assume $n\eq2$ in the introduction,} and $Z$ will be denoted by $C$. We assume moreover that all the singularities of the associated reduced curve $\Cr$ are {\it semi-weighted-homogeneous,} that is, there are analytic coordinates $u,v$ around any singular point $p_j$ of $\Cr$ such that a local defining function of $\Cr$ at $p_j$ is a $\C$-linear combination of monomials $u^{m_1}v^{m_2}$ with $w_jm_1\pl w'_jm_2\gess d_j$ and its lowest weighted-degree $d_j$ part, which is the $\C$-linear combination of $u^{m_1}v^{m_2}$ with $w_jm_1\pl w'_jm_2\eq d_j$, is {\it reduced\one} (since $\dim \Cr\eq1$). Here $w_j,w'_j$ are {\it mutually prime positive integers,} called the local (absolute) {\it weights,} and $d_j$ is a positive integer, called the (absolute) {\it weighted degree.} Let
$$f\eq\mprod_{k=1}^rf_k^{a_k}\in\C[x,y,z]$$
be a defining homogeneous polynomial of $C\sst Y$ with $f_k$ irreducible and multiplicity $a_k\gess1$ for $k\ins[1,r]$. Set\vskip-5mm
$$\aligned C_k\defs\{f_k\eq0\}\sst Y,&\q d'_k\defs\deg C_k\eq\deg f_k,\\ d\defs\deg C\eq\msum_{k=1}^r\,a_kd'_k,&\q d'\defs\deg\Cr\eq\msum_{k=1}^r\,d'_k.\endaligned$$
For $i\ins[1,d]$, put\vskip-5mm
$$\aligned\si_i&\defs\msum_{k=1}^r\,d'_k\bl(\lceil a_ki/d\rceil{-}1\br)\ins[0,i{-}1]\caps\Z,\\ \iota_i&\defs i-\si_i\ins[1,i]\caps\Z,\\ \be_{k,i}&\defs a_ki/d\mi\lceil a_ki/d\rceil\pl1\ins(0,1].\endaligned$$
\sk
Let $p_j$ ($j\ins[1,q]$) be the singular points of $\Cr$ which are {\it not ordinary double points.} Let $C_{j,l}$ ($l\ins[1,n_j]$) be the local irreducible components of $(\Cr,p_j)$ with $d_{j,l}\ins\Z$ the (absolute) weighted degree of $C_{j,l}$, which is either $w_j$ or $w'_j$ or $w_jw'_j$ with lowest weighted-degree part of the defining function respectively $u$ or $v$ or $u^{w'_j}\mi c_{j,l}y^{w_j}$ ($c_{j,l}\ins\C^*$), where $w_j,w'_j\ins\Z$ are the (absolute) weights of $\Cr$ at $p_j$, see Section\,\,\ref{S3} below. Let $a_{j,l}$ be the multiplicity of $(C,p_j)$ along $C_{j,l}$ (where $a_{j,l}\eq a_k$ if $C_{j,l}\sst C_k$). Set
$$\aligned\be_{j,l,i}&\defs a_{j,l}i/d\mi\lceil a_{j,l}i/d\rceil\pl1\in(0,1],\\ \ga_{j,i}&\defs\msum_{l=1}^{n_j}\,\be_{j,l,i}d_{j,l}\in(0,d_j],\\ N_j(n)&\defs\#\bl\{(m_1,m_2)\ins\Z_{>0}^2\mid w_jm_1\pl w'_jm_2\les n\br\}\,\,\,(n\ins\N).
\endaligned$$
Note that $d_j\eq\msum_l\,d_{j,l}$. Calculating the direct image by a weighted blowup, we can show the following.

\begin{ithm} \label{T2}
With the above notation and assumption, we have for integers $i\ins[1,d]$ the following  equalities\,$:$
\begin{equation} \label{4}
\aligned n_{f,\frac{i}{d}}&=\tbinom{\iota_i-1}{2}-\msum_{j=1}^q\,N_j\bl(\lceil\ga_{j,i}\rceil-1\br),\\ n_{f,\frac{i}{d}+1}&=(d'{-}3)d'\pl3-\msum_{j=1}^{q''}\bl(\tfrac{d_j}{w_j}{-}1\br)\bl(\tfrac{d_j}{w'_j}{-}1\br)\mi n_{f,\frac{i}{d}}\mi n_{f,\frac{i}{d}+2}\mi\de_{i,d},\\ n_{f,\frac{i}{d}+2}&=\tbinom{d'-\iota_i-1}{2}-\msum_{j=1}^q\,N_j\bl(d_j\mi\lceil\ga_{j,i}\rceil\br)-\de_{i,d}.\endaligned
\end{equation}
\end{ithm}

Here the $p_j$ ($j\ins[1,q'']$) are all the singular points of $\Cr$. Note that $n_{f,\al}\eq0$ unless $\al\ins(0,3)\caps\tfrac{1}{d}\one\Z$ (see \eqref{1.3} below), and $\binom{n}{2}\defs n(n{-}1)/2$ ($\forall\,n\ins\Z$). Theorem\,\,\ref{T2} says that if all the singularities of $\Cr$ are semi-weighted-homogeneous, the spectrum of the cone of a non-reduced plane curve $C$ can be expressed by using the local weights and the weighted degrees and multiplicities of local irreducible components together with the degrees and multiplicities of global ones without employing the relation between global components and singular points. This is a generalization of a formula for the reduced line arrangement case \cite{BS}, since\vskip-7mm
$$N_j(n)\eq\tbinom{n}{2}\q\h{for}\,\,\,n\ins\N,$$
in the {\it ordinary\one} (that is, semi-homogeneous) singularity case, where $w_j\eq w'_j\eq1$. In this case, the latter information is given by the $m_{j,k}\defs{\rm mult}_{p_j}C_k$, and the formula for $n_{f,\frac{i}{d}+e}$ with $e\eq0,2$ was essentially shown in \cite{Y} by computing the characteristic classes, since the $\ga_{j,i}$ can be expressed as $\msum_k\,m_{j,k}\be_{k,i}$. (The two formulas look apparently different, but it is not difficult to see that they are actually equivalent.) We can also give a simple proof of a formula which is equivalent to the one in \cite{Y} for $e\eq1$ by the following.

\begin{icor} \label{C2}
Assume all the singularities of $\Cr$ are ordinary. For $i\ins[1,d]\caps\Z$, we have
\begin{equation} \label{5}
\aligned n_{f,\frac{i}{d}+1}&=(\iota_i{-}1)(d'{-}\iota_i{-}1)+\msum_{k=1}^r\tbinom{d'_k}{2}\\&\q-\msum_{j=1}^q\bl(\lceil\ga_{j,i}\rceil{-}1\br)\bl(m_j{-}\lceil\ga_{j,i}\rceil\br)-\msum_{j=1}^{q'}\msum_{k=1}^r\tbinom{m_{j,k}}{2}.\endaligned
\end{equation}
\end{icor}

The last summation is taken over $k\ins[1,r]$ and $j\ins[1,q']$ with $\{p_j\}_{j\in[q+1,q']}$ the ordinary double points of $\Cr$ which are ordinary double point of some $C_k$, where $\binom{m_{j,k}}{2}\eq1$.

\begin{irem} \label{R1}
Assume all the irreducible components $C_k$ are smooth. Let $C'$ be the union of sufficiently small and general deformations $C'_k$ of $\,C_k$ $(k\ins[1,r])$ having only normal crossing singularities. Put $U\defs Y\stm\Cr$, $U'\defs Y\stm C'$. We can get
\begin{equation} \label{6}
\chi(U')-\tbinom{d'-2}{2}=\msum_{k=1}^r\tbinom{d'_k}{2},
\end{equation}
applying Remark\,\,\ref{R4} below to each $C'_k$, which is identified with a smooth deformation of a general line arrangement of degree $d'_k$, where $\nu_2\eq\binom{d'_k}{2}$. (One can prove \eqref{6} also by calculating $\chi(C')$.) Note that in the {\it line arrangement\one} case, we have
\begin{equation} \label{7}
\chi(U')\eq\tbinom{d'{-}2}{2}\,\bl({=}\,3\mi 2d'\pl\tbinom{d'}{2}\br).
\end{equation}
\end{irem}

\begin{irem} \label{R2}
The second equality of Theorem\,\,\ref{T2} follows from \eqref{2.3} and Remark\,\,\ref{R4} below. The second one of Corollary\,\,\ref{C1} is a consequence of the first and third using \eqref{2.3}, \eqref{6.1}, \eqref{6.2}, and Remark\,\,\ref{R4} below. Hence we do not have to study $\Om_{\Yt}^1(\log\Ctr)$, and the arguments are simple in the proof of Theorems\,\,\ref{T1} and \ref{T2}. Using the latter theorems and Corollary\,\,\ref{C2}, one can recover \cite[Theorem 3]{BS} for {\it reduced\one} line arrangements, where $d\eq d'$, $a_k\eq d'_k\eq1$, $\be_{k,i}\eq \be_{j,l,i}\eq i/d$, $\si_i\eq0$, $\iota_i\eq i$, $\ga_{j,i}\eq m_ji/d$ for any $i\ins[1,d]$, $k\ins[1,r]$, $j\ins[1,q]$.
\end{irem}

\begin{irem} \label{R3}
Theorems\,\,\ref{T1} and \ref{T2} imply that the coefficients $n_{f,\frac{i}{d}+e}$ for $i\ins[1,d]\caps\Z$, $e\eq0,2$ in the {\it reduced\one} case is determined only by the the spectral numbers $\al_{j,l}$ or the local weights $w_j,w'_j$ and degree $d_j$ for $j\ins[1,q]$ together with the global degree $d$, and is independent of $\nu_2$. (Here $\nu_m$ denotes the number of points of $\Cr$ with multiplicity $m\gess 2$.) Note that $\nu_2$ is reflected to the Euler number $\chi(U)$; for instance compare the generic line arrangement and smooth curve cases.
\end{irem}

\begin{irem} \label{R4}
Assume $C$ is {\it reduced\one} so that $d\eq d'$. Let $D\sst\PP^2$ be a smooth curve of degree $d$. We have the equality
\begin{equation} \label{8}
\chi(C)=\chi(D)+\msum_{j=1}^{q''}\one\mu_j,
\end{equation}
using a smoothing of $C$ in $Y\eq\PP^2$ defined by $h\defs f\pl sg$ ($s\ins\C$) together with a distinguished triangle $i_0^*\tos\psi_s\tos\varphi_s\,{\buildrel{[1]\,\,}\over{\to}}$. Here $g$ is a defining polynomial of $D$, and $i_0$ denotes the inclusion $\PP^2{\times}\{0\}\into\PP^2{\times}\C$. For the proof of \eqref{8} we {\it may assume\one} that $D$ intersects $C$ at {\it smooth\one} points, hence $\{h\eq0\}\sst\PP^2{\times}\C$ is {\it smooth\one} around $C\eq\{f\eq0\}\sst\PP^2{\times}\{0\}$. Note that
\begin{equation} \label{9}
\chi(D)\eq(3{-}d)d\,\,\bl({=}\,2\mi(d{-1})(d{-}2)\br).
\end{equation}
\end{irem}

\begin{irem} \label{R5}
The {\it pole order spectrum\one} $^P\Sp_f(t)\eq\msum_{\al}\,{}^P\!n_{f,\al}t^{\al}$ can be defined by replacing the Hodge filtration $F$ with the pole order filtration $P$ on the vanishing cohomology. We have the implications
\begin{equation} \label{10}
n_{f,\frac{3}{d}}\nes0\Longrightarrow{}^P\!n_{f,\frac{3}{d}}\nes0\Longrightarrow b_f\bl(-\tfrac{3}{d}\br)\eq0,
\end{equation}
where the converses do not necessarily hold, see \cite{bCM}, \cite{nwh}, \cite[Remark 4.2c]{wh}. One can easily verify that the {\it strong monodromy conjecture\one} for $f$ is reduced to the vanishing of $b_f\bl(-\frac{3}{d}\br)$ in the case the polynomial $f$ is reduced, essential, and indecomposable, and has only ordinary singularities, and moreover it follows from the $E_2$-degeneration of the pole order spectral sequence \cite{wh} in the {\it locally homogeneous\one} reduced curve case by showing the vanishing of $N_{d+3}$ as in \cite[Proposition 2]{wh}, see also \cite{BSY}, \cite{Wa} for the line arrangement case. Note that some examples with $b_f(-\tfrac{3}{d})\nes0$ are known in the {\it locally weighted homogeneous\one} reduced curve case, for instance $f\eq x^ay^{d-a}{+}z^d$ with $a\sgt d{-}a\sgt1$. There are many examples with $n_{f,\frac{3}{d}}\eq0$ even in the reduced curve case (see Sections\,\,\ref{S9}--\ref{S12} below), but no examples with $^P\!n_{f,\frac{3}{d}}\eq0$ are known in the {\it locally homogeneous\one} reduced curve case, see \cite[Example 5.5]{nwh} and Section\,\,\ref{S12} below for the non-reduced case. It is not necessarily easy to verify these, since a computation of pole order spectral sequence takes rather long.
\end{irem}
\sk
In Part~1 we review some basics of spectrum. In Part~2 we prove Theorems\,\,\ref{T1}--\ref{T2} and Corollary\,\,\ref{C2}. In Part~3 we explain some explicit computations of examples in the ordinary singularity case.

\tableofcontents
\numberwithin{equation}{section}

\part{Preliminaries} \label{Pa1}
\nin
In this part we review some basics of spectrum.

\section{Spectrum} \label{S1}
Let $f$ be a holomorphic function on a complex manifold $(X,0)$ with $f(0)\eq0$. Let $F_{\!f}$ be the Milnor fiber of $f$ at $0\ins X$. Let $i_0\col\{0\}\into X$ be the inclusion, The vanishing cohomology has a canonical mixed Hodge structure using the isomorphism{}
\begin{equation} \label{1.1}
H^{\ssb}(F_{\!f},\Q)=H^{\ssb}i_0^*\psi_f\Q_{h,X}.
\end{equation}
Here $\Q_{h,X}$ denotes the constant Hodge module of weight $d_X$ shifted by $-d_X$ whose underlying $\Q$-complex is the constant sheaf $\Q_X$, $\psi_f$ is the nearby cycle functor of mixed Hodge modules (preserving mixed Hodge modules up to shift by $1$), and $H^{\ssb}i_0^*$ is the {\it cohomological\one} pullback functor of mixed Hodge modules, see \cite{mhm}.
\sk
We define the {\it spectrum\one} $\Spf(t)=\sum_{\al\in\Q}n_{f,\al}\one t^{\al}$ by
\begin{equation} \label{1.2}
\aligned&n_{f,\al}:=\msum_{j=0}^{d_X-1}(-1)^j\dim_{\C}\Gr_F^p\Ht^{d_X-j-1}(F_{\!f},\C)_{\la}\\&\q\q\h{with}\q p\eq[d_X\mi\al],\,\la=e^{-2\pi\sqrt{-1}\al},\endaligned
\end{equation}
as in \cite{BS0}, \cite{BS}. Here $\Ht^{\ssb}$ is the reduced cohomology, and the lower index $_{\la}$ indicates the $\la$-eigenspace of the action of $T_s$ with $T_s$ the semisimple part of the Jordan decomposition $T\eq T_sT_u$ of the monodromy $T$. In the isolated singularity case, the {\it spectral numbers\one} $\al_{f,1},\dots,\al_{f,\mu_f}$ (with $\mu_f$ the Milnor number) are defined to be the rational numbers $\al$ with $n_{f,\al}\nes0$ (counted with multiplicity).

\begin{rem} \label{R1.1}
The above definition of spectrum and spectral numbers is compatible with \cite{St0} in the isolated singularity case, although it is shifted by 1 compared with the one in \cite{St}. It is known (see \cite[Proposition 5.2]{BS0}) that with our definition we have
\begin{equation} \label{1.3}
n_{f,\al}=0\q\h{if}\,\,\,\al\notin(0,d_X).
\end{equation}
\sk
The spectrum depends only on the hypersurface $f^{-1}(0)$, since so does the embedded resolution used for the calculation of the mixed Hodge structure.
\end{rem}

\begin{rem} \label{R1.2}
There is a ``dual" definition of spectrum, which is denoted by $\Spf'(t)$ (see \cite{ste}), and we can show the equality
\begin{equation} \label{1.4}
\Spf'(t)=\Spf(t^{-1})\one t^{d_X}.
\end{equation}
In the isolated singularity, we have 
\begin{equation} \label{1.5}
\Spf'(t)=\Spf(t),
\end{equation}
using the self-duality of the vanishing cohomology. This is equivalent to the ``symmetry" of spectral numbers $\al_{f,1},\dots,\al_{f,\mu_f}$ if these are indexed weakly increasingly.
\end{rem}

\begin{rem} \label{R1.3}
For $f\ins\C\{x\}$, $g\ins\C\{y\}$ having isolated singularities at 0 with $f(0)\eq g(0)\eq0$, we have the Thom--Sebastiani-type theorem
\begin{equation} \label{1.6}
\Sp_{\mhs f+g}(t)=\Spf(t)\one {\cdot}\one\Sp\mhs_g(t),
\end{equation}
meaning that the spectral numbers of $f\pl g\ins\C\{x,y\}$ are the sums of spectral numbers of $f$ and $g$, that is, $\al_{f,i}\pl\al_{g,j}$ for $i\ins[1,\mu_f],j\ins[1,\mu_g]$, where $\mu_{f+g}\eq\mu_f\one\mu_g$, see \cite{SS}, \cite{Va0}.
\end{rem}

\section{Homogeneous case} \label{S2}
Assume $f$ is a homogeneous polynomial with $X\eq\C^{d_X}$. Let $\rho\col\Xt\tos X$ be the blowup at 0 with exceptional divisor $E\cong\PP^{d_X-1}$. Note that $\Xt$ is a line bundle over $E$. Let $D\sst E$ be the intersection of $E$ and the proper transform of $f^{-1}(0)$, which is identified with the projective hypersurface defined by $f$. Let $\pi\col(\Et,\Dt)\tos(E,D)$ be an embedded resolution of $D\sst E$. We get an embedded resolution $\rhot$ of $f^{-1}(0)\sst X$ composing $\rho\col\Xt\tos X$ with the base change $\pit$ of $\pi\col\Et\tos E$ by the canonical projection of the {\it line bundle\one} $\Xt\onto E$ as in the diagram bellow:
$$\xymatrix{X & \Xt \ar[l]_{\rho} \ar[d] & \Xt{\times}_E\Et \ar[l]_{\pit} \ar[d]\\
& E & \Et \ar[l]_{\pi}}$$
\sk
Set $U\defs E\stm D=\Et\stm\Dt$ with $\jt_U\col U\into\Et$ and $\iti_{\Et}\col\Et\into\Xt{\times}_E\Et$ natural inclusions, where the latter is the zero section of line bundle. Let $\ft$ be the pullback of $f$ on $\Xt_{\Et}\defs\Xt{\times}_E\Et$. From \cite[Theorem 4.2]{BS0} we can deduce the isomorphism
\begin{equation} \label{2.1}
\iti_{\Et}^*(\psi_{\ft}\Q_{h,\Xt_{\Et}}[d_X{-}1])=\jt_{U*}\jt_U^*\iti_{\Et}^*(\psi_{\ft}\Q_{h,\Xt_{\Et}}[d_X{-}1]),
\end{equation}
where the right-hand side is a mixed Hodge module, since $\jt_U$ is an affine morphism. We have moreover the following isomorphisms in the bounded derived category of {\it algebraic\one} mixed Hodge modules (see \cite[\S4]{mhm}):
\begin{equation} \label{2.2}
\aligned&\rhot'_*\iti_{\Et}^*(\psi_{\ft}\Q_{h,\Xt_{\Et}}[d_X{-}1])=i_0^*\rhot_*(\psi_{\ft}\Q_{h,\Xt_{\Et}}[d_X{-}1])\\&=i_0^*(\psi_f\rhot_*\Q_{h,\Xt_{\Et}}[d_X{-}1])=i_0^*(\psi_f\Q_{h,X}[d_X{-}1]),\endaligned
\end{equation}
where $\rhot'$ is the restriction of $\rhot$ to $\Et\eq\rhot^{-1}(0)$. For the last isomorphism, we apply the decomposition theorem saying that $\rhot_*\Q_{h,\Xt_{\Et}}[d_X]$ is isomorphic to the direct sum of $\Q_{h,X}[d_X]$ and a complex of mixed Hodge modules supported on $f^{-1}(0)$. Note that
$$\jt_U^*\iti_{\Et}^*(\psi_{\ft}\Q_{h,\Xt_{\Et}}[d_X{-}1])$$
is a locally constant variation of Hodge structures of type $(0,0)$ and rank $d\defs\deg f$ on $U$. The monodromy eigenvalues are $d\one $th root of unity, and the eigensheaves with complex coefficients have rank one. So we get local systems $L_i$ of rank one with complex coefficients for $i\ins[1,d]\caps\Z$ such that
\begin{equation} \label{2.3}
H^{\ssb}(U,L_i)=H^{\ssb}(F_{\!f},\C)_{\la}\q\h{with}\q\la\defs e^{-2\pi\sqrt{-1}\,i/d},
\end{equation}
see also \cite{Di}, \cite[\S1.6]{BS} (where the index is given in a different way). This is used to show the middle equalities of \eqref{3} and \eqref{4}. The {\it Hodge filtration\one} $F$ on $H^{\ssb}(U,L_i)$ can be calculated by using the truncations $\sigma_{\ges\ssb}$ of the logarithmic complexes on $\Et$ (applying \cite[(3.11.2)]{mhm} inductively). Note that the monodromy $T$ is semisimple, since the {\it geometric monodromy\one} $\ga^*$ is the inverse of the monodromy $T$, where $\ga$ is defined by $(x_0,\dots,x_n)\mapsto(\zeta_dx_0,\dots,\zeta_dx_n)$ with $\zeta_d\defs e^{2\pi\sqrt{-1}/d}$, see \cite{DiSa0}.
\sk
The {\it monodromy\one} of the local system $L_i$ around a general point of an irreducible component $\Dt_k$ of $\Dt$ is given by multiplication by $e^{2\pi\sqrt{-1}\,a_ki/d}$ with $a_k$ the multiplicity of $\Dt$ along $\Dt_k$. This can be verified, since $\ft$ is the product of $t_0^d$ and $t_k^{a_k}$ at a sufficiently general point of $\Dt_k$, where $t_0$ and $t_k$ are local defining functions of $\Et$ and the inverse image of $\Dt_k$ by the projection $\Xt_{\Et}\to\Et$ respectively, see also \cite[Theorem 3.3]{mhm}.

\section{Semi-weighted-homogeneous case} \label{S3}
Let $h$ be a germ of a holomorphic function of $n$ variables having an isolated singularity. Assume $h$ is {\it semi-weighted-homogeneous\one} so that its (absolute) {\it weights\one} are $w_i\ins\Z_{>0}$ and its {\it lowest weighted degree\one} is $d\sgt1$. Here the greatest common divisor of the $w_i$ is one, and the lowest weighted-degree $d$ part of $h$ has an isolated singularity, see the introduction for the two-variable case. It is well known that the spectrum of $h$ is expressed as
\begin{equation} \label{3.1}
\Sp_h(t)=\mprod_{i=1}^n\bl(t\mi t^{w_i/d}\br)\big/\bl(t^{w_i/d}\mi1\br).
\end{equation}
This can be reduced to the weighted homogeneous case in \cite{St0} using \cite{Va} or \cite{DMST} together with the $\mu$-constant deformation
\begin{equation} \label{3.2}
h_s\defs\msum_{m\ges d}\,s^{m-d}h^{(m)}\q(s\ins S).
\end{equation}
Here $S\sst\C$ is an open disk with $1\ins S$, and $h\eq\msum_{m\ges d}\,h^{(m)}$ with $h^{(m)}$ weighted homogeneous of weighted degree $m\gess d$, that is,
\begin{equation} \label{3.3}
h^{(m)}(s^{w_1}z_1,\dots,s^{w_n}z_n)\eq s^mh^{(m)}(z_1,\dots,z_n).
\end{equation}
We may assume that $h$ is a polynomial by the finite determinacy, see \cite{GLS}. One can see that it is a $\mu$-constant deformation by taking the graded quotients of the Koszul complex $(\Om_{\C^n,0}^{\ssb},{\rm d}h_s\wedge)$ by the decreasing filtration defined by weighted degrees and using the Mittag-Leffler condition together with the exactness of the completion functor of finite modules over the convergent power series.
\sk
Setting $t'\defs t^{1/d}$, $t''\defs t'\mi 1$, and taking the limit for $t''\to0$ in \eqref{3.1}, we get that
\begin{equation} \label{3.4}
\mu_h=\mprod_{i=1}^n(d\mi w_i)/w_i.
\end{equation}

\begin{rem} \label{R3.1}
Assume $n\eq2$. Each irreducible factor $h_l$ of $h$ is semi-weighted-homogeneous with (absolute) weights $w_1,w_2$. Indeed, the product of the lowest graded-degree $d_l$ part $h'_l$ of $h_l$ coincides with that of $h$, which is reduced by hypothesis.
\sk
Assume further that $h$ is {\it irreducible.} Then the weighted-degree $d$ part $h'$ of $h$ can be written as $z_1^{w_2}\pl z_2^{w_1}$ for some local analytic coordinates $z_1,z_2$. Indeed, $h'$ is a linear combination of monomials with weighted degree $d$ by hypothesis. Moreover, it is not divisible by $z_1$ nor $z_2$ (since $d\sgt1$). This can be shown for instance using a $\mu$-constant deformation or a formula for spectrum in the Newton-nondegenerate two-variable case, see \cite{St0}. (Recall that the number of local irreducible components minus one coincides with the dimension of the unipotent monodromy part of the vanishing cohomology.)
Take now the pull-back of $h$ by the ramified covering defined by $(z_1,z_2)\mapsto(z_1^{w_1},z_2^{w_2})$. Its lowest (homogeneous) degree part, which coincides with the pullback of $h'$, is {\it homogeneous\one} and {\it reduced\one} (hence it is a product of mutually distinct linear functions), since the map is locally biholomorphic outside $\{z_1z_2\eq1\}$. This implies the desired assertion, since the degree of the ramified covering is $w_1w_2$, the image of each irreducible component is irreducible, and $\{h\eq0\}$ is irreducible.
\end{rem}

\part{Proof of the main theorems} \label{Pa2}
\nin
In this part we prove Theorems\,\,\ref{T1}--\ref{T2} and Corollary\,\,\ref{C2}.

\section{Proof of the first equality in Theorem~\ref{T2}} \label{S4}
Let $\Lc_i$ be the Deligne extension of $\OO_U{\otimes}_{\C}L_i$. More precisely, the restriction of $\Lc_i$ to the open subset $U''\defs Y\stm\mcup_{j=1}^q\{p_j\}$ is the Deligne extension of $\OO_U{\otimes}_{\C}L_i$ over $U''$ such that the residues of connection are contained in $[0,1)$ (see \cite{De}), and $\Lc_i$ is its direct image by the inclusion $j''\col U''\into Y$. This is a locally free sheaf of rank 1, and is locally isomorphic to the tensor product of the Deligne extensions $\Lc_{j,l,i}$ of local systems $L_{j,l,i}$ defined on the complements of local irreducible components $C_{j,l}$ of $\Cr$ at each $p_j$. Note that each $\Lc_{j,l,i}$ is generated by the product of $h_{j,l}^{\al_{j,l,i}}$ with a locally constant section, where $h_{j,l}$ is a local defining function of $C_{j,l}$ and $\al_{j,l,i}\eq1{-}\be_{j,l,i}$. (In this paper we use analytic sheaves, and $Y$ is a complex manifold.)
\sk
Since $1\mi\be_{k,i}\ins[0,1)$ is the residue of the logarithmic connection at general points of $C_k$, we see that
\begin{equation} \label{4.1}
\deg\Lc_i=-\msum_{k=1}^r\,d'_k(1\mi\be_{k,i})=-d'\pl\iota_i\,\in\,(-d',0],
\end{equation}
considering the $d$-fold self-tensor product of $\Lc_i$ if necessary. (Note that the degree of a line bundle on $\PP^2$ can be computed by restricting it to a general line in $\PP^2$.) We then get
\begin{equation} \label{4.2}
\deg\Om_Y^2(\log\Cr){\otimes}_{\OO_Y}\Lc_i=\iota_i\mi3\,\ges\,-2,
\end{equation}
where $\Om_Y^2(\log\Cr)\defs\Om_Y^2(\Cr)$, hence
\begin{equation} \label{4.3}
\chi\bl(\Om_Y^2(\log\Cr){\otimes}_{\OO_Y}\Lc_i\br)=\tbinom{\iota_i\mi1}{2}.
\end{equation}
Note that $\chi\bl(\OO_{\PP^2}(m)\br)=\binom{m+2}{2}$ for $m\in\Z$.
\sk
Let $w_j,w'_j\ins\Z_{>0}$ be the (absolute) weights of a semi-weighted-homogeneous singularity of $\Cr\sst Y$ at $p_j$. There is a local defining function $h_j$ of $\Cr$ which is a $\C$-linear combination of monomials $u^{m_1}v^{m_2}$ with $m_1w_j\pl m_2w'_j\ges d_j$ and such that its weighted-degree $d_j$ part $h'_j$, which is by definition the linear combination of monomials with $m_1w_j\pl m_2w'_j\eq d_j$, has an isolated singularity, where $u,v$ are local analytic coordinates of $Y^{\rm an}$ around $p_j$. We have the factorization $h_j\eq\mprod_l\one h_{j,l}$ such that $h_{j,l}$ is a local defining function of $C_{j,l}$ and it is semi-weighted-homogeneous of weighted degree$\,\gess d_{j,l}$ with $\msum_l\one d_{j,l}\eq d_j$, where we have the equality $h'_j\eq\mprod_l\one h'_{j,l}$ with $h'_{j,l}$ the lowest weighted-degree $d_{j,l}$-part of $h_{j,l}$. Note that $h'_j$ is {\rm reduced,} see Remark\,\,\ref{R3.1}.
\sk
We have a local $\C^*$-action defined by $(u,v)\mapsto(\la^{w_j}u,\la^{w'_j}v)$ for $\la\ins\C^*$ with $|\la|\less1$. Let $\pi\col\Yt\tos Y$ be the weighted blowup of $Y$ along $p_j$ for the weights $w_j,w'_j$. This can be obtained locally on $Y$ by taking the quotient by the covering transformation group of the usual blowup of the finite ramified covering in Remark\,\,\ref{R3.1}. Let $E\defs\pi^{-1}(p_j)$ be the exceptional divisor. (Here the index $_j$ is omitted to simplify the notation.) This is identified with the local $\C^*$-orbits, which are defined by $\{ru^{w'_j}\eq sv^{w_j}\}$ for $r,s\ins\C$ and is isomorphic to $\PP^1$. Note that the pullback of a meromorphic function $u^{m_1}v^{m_2}$ with $m_1,m_2\ins\Z$ has zeros of order 1 along the exceptional divisor $E$ in the case $m_1w_j\pl m_2w'_j\eq1$ (using $(w_j,w'_j)\eq1$).
\sk
The above construction gives an embedded resolution of $(Y,\Cr)$ by a $V$-manifold with $V$-normal crossings. Set $\Ct\defs\pi^{-1}(C)$. Let $\Lct_i$ be the Deligne extension of $\OO_U{\otimes}_{\C}L_i$ over $\Yt$. Since $\al_{j,l,i}\eq1\mi\be_{j,l,i}$ and
$$\msum_{l=1}^{n_j}\,(1\mi\be_{j,l,i})d_{j,l}\eq d_j\mi\ga_{j,i}\ges0,$$
we can verify on a neighborhood of $E$ that
\begin{equation} \label{4.4}
\aligned&\Lct_i=\OO_{\Yt}(eE){\otimes}_{\OO_{\Yt}}\pi^*\Lc_i\q\,\,\,\,\h{with}\\&e\defs\lfloor d_j\mi\ga_{j,i}\rfloor\eq d_j\mi\lceil\ga_{j,i}\rceil\ges0,\endaligned
\end{equation}
considering the pull-back of a generator of $\Lc_i$ at $p_j$, which is obtained by using the tensor product of Deligne extensions as explained above, see also \cite[\S3.2]{BaSa}. Let $\Ct'$ be the intersection of $E$ with the closure of $\Ctr\stm E$.
\sk
Assume for the moment the singularity is {\it ordinary,} that is, $w_j\eq w'_j\eq1$. Taking the {\it residue\one} along $E\sst\Yt$, we then see that
\begin{equation} \label{4.5}
\aligned\bl(\Om_{\Yt}^2(\log\Ctr){\otimes}_{\OO_{\Yt}}\Lct_i\br)|_E&\cong\Om_E^1(\log\Ct'){\otimes}_{\OO_E}\OO_{\Yt}(eE)|_E\\&\cong\OO_{\PP^1}\bl(\lceil\ga_{j,i}\rceil\mi2\br),\endaligned
\end{equation}
since $\deg\Om_E^1(\log\Ct')\eq d_j\mi2$ and $\OO_{\Yt}(E)|_E\cong\OO_{\PP^1}(-1)$. Hence
\begin{equation} \label{4.6}
\aligned\bl(\Om_{\Yt}^2(\log\Ctr)&{\otimes}_{\OO_{\Yt}}\Lct_i\br){\otimes}_{\OO_{\Yt}}\bl(\OO_{\Yt}(pE)/\OO_{\Yt}\bl((p{-}1)E\br)\br)\\&\cong\OO_{\PP^1}\bl(\lceil\ga_{j,i}\rceil\mi p\mi2\br)\q\q(\forall\,p\in\Z).\endaligned
\end{equation}
Here $\OO_{\Yt}(-pE)\eq\I_E^p$ for $p\ins\N$ with $\I_E\sst\OO_{\Yt}$ the ideal of $E\sst\Yt$. We then get
\begin{equation} \label{4.7}
R^n\pi_*\bl(\Om_{\Yt}^2(\log\Ctr){\otimes}_{\OO_{\Yt}}\Lct_i\br)\eq0\q\q(\forall\,n>0),
\end{equation}
restricting to $p\less0$ in \eqref{4.6}, see \cite[Ch.\,III, Theorem 11.1]{Ha}.
\sk
We then recognize that the vanishing \eqref{4.7} holds also in the case $(\Cr,p_j)$ is {\it weighted-homogeneous,} taking the invariant part of the action of the covering transformation group $G\defs(\Z/w_j\Z)\one {\times}\one(\Z/w'_j\Z)$ of the finite ramified covering $\rho\col Y'\to Y$ explained in Remark\,\,\ref{R3.1} (where $w_1\eq w_j$, $w_2\eq w_j'$). This group acts on the direct images by $\rho$ and $\rhot$ in the following commutative diagram
$$\xymatrix{\Yt \ar[d]^{\pi} & \Yt' \ar[l]_{\rhot} \ar[d]^{\pi'}\\
Y & Y' \ar[l]_{\rho}}$$
Here $\pi'$ is the usual blowup, and $Y$ denotes a sufficiently small polydisk with center $p_j$. (The latter may be replaced by an affine space for the construction of the quotient of the usual blowup by the action of $G$.)
\sk
To deduce the vanishing \eqref{4.7} for the {\it semi-weighted-homogeneous\one} case, we consider the $\mu$-constant deformation $h_s$ ($s\ins S$) as in Section\,\,\ref{S3}. We have the reduced divisor $C_S\sst\YS\defs Y{\times}S$ defined by the $h_s$. Here $Y$ denotes a small  polydisk with center $p_j$ as above, and $h_j$ is denoted by $h$ to simplify the notation. Let $\Ct_S\sst\YtS\defs\Yt{\times}S$ be the proper transform of $C_S\sst\YS$ by the morphism $\pit\defs\pi{\times}{\rm id}\col\YtS\to\YS$. The Deligne extensions $\Lc_i$, $\Lct_i$ on $Y{\times}\{1\}$, $\Yt{\times}\{1\}$ can be extended over $\YS$ and $\YtS$, which are denoted by $\Lc_{S,i}$ and $\Lct_{S,i}$ respectively. For $c\ins S$, set $C_c\defs C_S\cap(Y{\times}\{c\})$, $\Ct_c\defs\Ct_S\cap(\Yt{\times}\{c\})$. Let $\Lc_{c,i}$, $\Lct_{c,i}$  be the restrictions of $\Lc_{S,i}$, $\Lct_{S,i}$ to $Y{\times}\{c\}$, $\Yt{\times}\{c\}$ respectively. Put 
$$\aligned\Om^p_{\rm rel,log}(\Lc_{S,i})&\defs\Om^p_{\YS/S}(\log C_S){\otimes}_{\OO}\Lc_{S,i},\\ \Om^p_{\log}(\Lc_{c,i})&\defs\Om^p_Y(\log C_c){\otimes}_{\OO}\Lc_{c,i},\endaligned$$
and similarly for $\Om^p_{\rm rel,log}(\Lct_{S,i})$, $\Om^p_{\log}(\Lct_{c,i})$, replacing $\Lc,Y,C$ with $\Lct,\Yt,\Ct$. We have the long exact sequence
\begin{equation} \label{4.8}
\to R^n\pit_*\Om^2_{\rm rel,log}(\Lct_{S,i})\buildrel{\!\!s-c}\over\longrightarrow R^n\pit_*\Om^2_{\rm rel,log}(\Lct_{S,i})\to R^n\pi_*\Om^2_{\log}(\Lct_{c,i})\to
\end{equation}
where $s$ denotes the coordinate of $S\sst\C$. Using Nakayama's lemma, we see that the vanishing \eqref{4.7} for $R^n\pi_*\Om^2_{\log}(\Lct_{c,i})$ ($n\sgt0$) with $c\eq0$ implies that for $R^n\pit_*\Om^2_{\rm rel,log}(\Lct_{S,i})$ ($n\sgt0$) with $S$ sufficiently small, and then for $R^n\pi_*\Om^2_{\log}(\Lct_{c,i})$ ($n\sgt0$) with $c\ins S$. Using the local $\C^*$-action as in \eqref{3.3}, the last vanishing holds also for $c\eq1$. So the vanishing \eqref{4.7} is proved in the semi-weighted-homogeneous case.
\sk
The above argument implies for $c\ins S$ the short exact sequence
\begin{equation} \label{4.9}
0\to\pit_*\Om^2_{\rm rel,log}(\Lct_{S,i})\buildrel{\!\!s-c}\over\longrightarrow\pit_*\Om^2_{\rm rel,log}(\Lct_{S,i})\to\pi_*\Om^2_{\log}(\Lct_{c,i})\to0.
\end{equation}
Applying the snake lemma to the injective morphism to the short exact sequence
\begin{equation} \label{4.10}
0\to\Om^2_{\rm rel,log}(\Lc_{S,i})\buildrel{\!\!s-c}\over\longrightarrow\Om^2_{\rm rel,log}(\Lc_{S,i})\to\Om^2_{\log}(\Lc_{c,i})\to0,
\end{equation}
we see that the cokernel of the inclusion
$$\pit_*\Om^2_{\rm rel,log}(\Lct_{S,i})\into\Om^2_{\rm rel,log}(\Lc_{S,i}$$
is locally free, and that of
$$\ep_{c,i}\col\pi_*\Om^2_{\log}(\Lct_{c,i})\into\Om^2_{\log}(\Lc_{c,i}),$$
which is identified with a vector space $E_{c,i}$, is independent of $c\ins S$, using the local $\C^*$-action as in \eqref{3.3}.
\sk
By \eqref{4.3} and \eqref{4.7}, the first equality of Theorem\,\,\ref{T2} is then reduced to
\begin{equation} \label{4.11}
\dim E_{0,i}=N_j\bl(\lceil\ga_{j,i}\rceil\mi1\br).
\end{equation}
Since $\Om^2_Y(\log C_0)\eq\Om^2_Y(C_0)$ with weighted degree of the defining function $h'$ of $C_0$ equal to $d_j$ and the pullback of the ``meromorphic" form $\tfrac{\ddd u}{u}{\wedge}\tfrac{\ddd v}{v}$ generates the logarithmic 2-forms {\it at general points\one} of the exceptional divisor $E$, the assertion follows using \eqref{4.4}. Indeed, the product of the local generator of $\Lc_i$ at $p_j$ as above and a logarithmic 2-form
$$h'{}^{-1}u^{m_1-1}v^{m_2-1}\ddd u{\wedge}\ddd v$$
belongs to the image of the injection $\ep_{0,i}$ if and only if
\begin{equation} \label{4.12}
(d_j\mi\ga_{j,i})\mi d_j\pl w_jm_1\pl w'_jm_2\gess0,
\end{equation}
calculating the vanishing order restricted to a union of general local $\C^*$-orbits in the blowup. Note that the {\it Lie derivative\one} by the vector field corresponding to the local $\C^*$-action can be computed by using the integral curves of the vector field, that is, the orbits of the local $\C^*$-action, see also \cite{BaSa}. This finishes the proof of the first equality in Theorem\,\,\ref{T2}.

\section{Proof of the third equality in Theorem~\ref{T2}} \label{S5}
Let $L^{\vee}_i$ be the dual local system of $L_i$, and $\Lc_i^{\vee}$, $\Lct_i^{\vee}$ be its Deligne extension over $\Yt$ and the direct image by $j''$ of its Deligne extension over $U''$ respectively, where the residues of connection are contained in $(0,1]$ instead of $[0,1)$ so that there are isomorphisms of locally free sheaves
\begin{equation} \label{5.1}
\Lct^{\vee}_i(\Ctr)\eq\Hc om_{\OO}(\Lct_i,\OO_{\Yt}),\q\Lc^{\vee}_i(\Cr)\eq\Hc om_{\OO}(\Lc_i,\OO_Y),
\end{equation}
using the Hartogs extension theorem for the last isomorphism (since the latter holds on $U''$). Here $(\Ctr)$ means the tensor product with $\OO_{\Yt}(\Ctr)$. Moreover, the filtered logarithmic complexes
$$\bl(\Om_{\Yt}^{\ssb}(\log\Ctr){\otimes}_{\OO}\Lct^{\vee}_i,F\br),\q\bl(\Om_{\Yt}^{\ssb}(\log\Ctr){\otimes}_{\OO}\Lct_i,F\br)$$
are the duals (defined by using the dualizing sheaf $\Om^2_{\Yt}$) of each other up to shifts of degrees of filtration and complex, see \cite[Section 2.4]{mhp}. Here we can apply also \cite[Proposition 3.11]{mhm}. Anyway we get
\begin{equation} \label{5.2}
\aligned\chi(\Lc_i)&=\chi\bl(\Om_Y^2(\log\Cr){\otimes}_{\OO_Y}\Lc^{\vee}_i\br),\\\chi(\Lct_i)&=\chi\bl(\Om_{\Yt}^2(\log\Ctr){\otimes}_{\OO_{\Yt}}\Lct^{\vee}_i\br),\endaligned
\end{equation}
where we have to use \eqref{5.1} for the first equality.
\sk
By the argument in Section\,\,\ref{S4} the residue of connection of $\Lc_i$ at general points of $C_k$ is given by $\be_{k,i}\ins(0,1]$, hence we get
\begin{equation} \label{5.3}
\deg\Lc^{\vee}_i=-\msum_{k=1}^r\,d'_k\be_{k,i}=-\iota_i,
\end{equation}
and
\begin{equation} \label{5.4}
\aligned\deg\Om_Y^2(\log\Cr){\otimes}_{\OO_Y}\Lc^{\vee}_i=d'\mi\iota_i\mi3,\\ \chi\bl(\Om_Y^2(\log\Cr){\otimes}_{\OO_Y}\Lc^{\vee}_i\br)=\tbinom{d'-\iota_i\mi1}{2}.\endaligned
\end{equation}
\sk
Corresponding to \eqref{4.4}, we have
\begin{equation} \label{5.5}
\aligned&\Lct^{\vee}_i=\OO_{\Yt}(eE){\otimes}_{\OO_{\Yt}}\pi^*\Lc^{\vee}_i\\ &\h{with}\q e\defs\lceil\ga_{j,i}\rceil\mi1\ges0.\endaligned
\end{equation}
Note that $\lceil\ga_{j,i}\rceil\mi1$ is the largest integer {\it strictly smaller than\one} $\ga_{j,i}\defs\msum_l\one\be_{j,l,i}d_{j,l}$, and hence $\ga_{j,i}\mi(\lceil\ga_{j,i}\rceil\mi1)\ins(0,1]$. We then get the vanishing of the higher direct images corresponding to \eqref{4.7} in the semi-weighted-homogeneous case.
\sk
As for \eqref{4.12}, we obtain the corresponding condition
\begin{equation} \label{5.6}
\ga_{j,i}\mi d_j\pl w_jm_1\pl w'_jm_2\sgt0,
\end{equation}
where the strict inequality holds, since the residue of connection at general points of $E$ belongs to $(0,1]$ instead of $[0,1)$. Since the condition \eqref{5.6} is equivalent to
\begin{equation} \label{5.7}
w_jm_1\pl w'_jm_2\gess d_j\mi\lceil\ga_{j,i}\rceil\pl1,
\end{equation}
the third equality in Theorem\,\,\ref{T2} follows.
\sk
This completes the proof of Theorem\,\,\ref{T2}, since it is sufficient to prove the first and third equalities by Remarks\,\,\ref{R2} and \ref{R4}.

\section{Proof of Corollary~\ref{C2}} \label{S6}
It is enough to show that $n_{f,\frac{i}{d}+1}$ in \eqref{4} coincides with the one in \eqref{5}. The summations over $j$ in the formulas of Theorem\,\,\ref{T2} and Corollary\,\,\ref{C2} can be extended to the summations over all the singular points of $\Cr$ including all the ordinary double points, where the corresponding summands vanish.
\sk
We have the identity
\begin{equation} \label{6.1}
\tbinom{x}{2}+\tbinom{y}{2}+xy=\tbinom{x+y}{2},
\end{equation}
which implies inductively that $\msum_k\tbinom{x_k}{2}\pl\msum_{k<k'}x_kx_{k'}\eq\tbinom{\sum_kx_k}{2}$, where $x_k\eq d'_k$ or $m_{j,k}$.
\sk
Applying \eqref{6.1} also to the case $x\eq\iota_i\mi1$, $y\eq d'\mi\iota_i\mi1$ or $x\eq\lceil\ga_{j,i}\rceil\mi1$, $y\eq m_j\mi\lceil\ga_{j,i}\rceil$, the assertion is then reduced to the equalities
\begin{equation} \label{6.2}
\tbinom{d'}{2}\pl\tbinom{d'-2}{2}\eq d'{}^2{-}3d'{+}3,\q\tbinom{m_j}{2}\pl\tbinom{m_j-1}{2}\eq(m_j{-}1)^2.
\end{equation}
Indeed, we have $m_j\eq\msum_k\,m_{j,k}$, and the intersection number $C_k\one {\cdot}\one C_{k'}$ can be calculated {\it globally\one} and {\it locally,} that is,
\begin{equation} \label{6.3}
C_k\one {\cdot}\one C_{k'}=d'_kd'_{k'}=\msum_{j=1}^{q''}\,m_{j,k}\one m_{j,k'}\q\h{if}\,\,\,k\nes k',
\end{equation}
where $\{p_j\}_{j\in[1,q'']}$ is the set of singular points of $\Cr$, see for instance \cite[\S8.4]{Fu}. So the assertion follows. This finishes the proof of Corollary\,\,\ref{C2}.

\begin{rem} \label{R6.1}
In the case $\Cr$ has only ordinary singularities and all the global irreducible components $C_k$ are smooth, we have a geometric proof of Corollary\,\,\ref{C2} as follows. Since the deformation is small and general, and $\Cr$ has multiplicity $m_j$ at $p_j$, this point becomes $m_j(m_j{-}1)/2$ ordinary double points under the deformation in Remark\,\,\ref{R1}. We then get
\begin{equation} \label{6.4}
\aligned\chi(\Cr)-\chi(C')&=\msum_{j=1}^q\,\bl(1\mi\bl(m_j\mi m_j(m_j{-}1)/2\br)\br)\\&=\msum_{j=1}^q\,\,(m_j{-}1)(m_j{-}2)/2.\endaligned
\end{equation}
Define $n'_{f,\frac{i}{d}+e}$ for $e\eq1$ and $e\eq0,2$ by the right-hand side of \eqref{5} and \eqref{4} respectively for the same index ${\frac{i}{d}\pl e}$. Set
$$\chi(U)'_i\defs\msum_{e=0}^2\,n'_{f,\frac{i}{d}+e}+\de_{i,d}.$$
By Theorem\,\,\ref{T2} we have to show the equality
\begin{equation} \label{6.5}
\chi(U)'_i=\chi(U)\q(\forall\,i\ins[1,d]\caps\Z).
\end{equation}
\sk
Using \eqref{6.1} for $x\eq\lceil\ga_{j,i}\rceil\mi1$, $y\eq m_j\mi\lceil\ga_{j,i}\rceil$ together with \eqref{6.4}, the assertion is reduced to the case $\Cr\eq C'$ (replacing $\Cr$ with $C'$), where $q\eq0$ and $\nu_m\eq0$ for $m\gess 3$. By \eqref{6.1} for $x\eq\iota_i\mi1$, $y\eq d'-\iota_i\mi1$, the assertion now follows from Remark\,\,\ref{R1}.
\end{rem}

\section{Proof of the first part of Theorem~\ref{T1}} \label{S7}
We may assume that $Z$ is reduced so that $Z\eq Z'$, $f\eq f'$, $d\eq d'$. If $Z$ is smooth (where $M_j(\tfrac{i}{d})\eq0$ for any $i,j$), the assertion follows from \eqref{3.1}. So it is enough to calculate the ``contribution of singular points".
\sk
Let $Z''\sst Y\eq\PP^n$ be a smooth hypersurface of degree $d$ defined by a polynomial $g$ and intersecting $Z$ at smooth points as in Remark\,\,\ref{R4}. We have a one-parameter deformation $f\!_S\defs f\pl sg$ defined on $\YS\defs Y{\times}S$, giving a smoothing $Z_c\defs Z_S\cap(Y{\times}\{c\})$ ($c\ins S$) with $Z_S\defs\{f\!_S\eq0\}\sst\YS$, where $S$ is a Zariski-open neighborhood of $0\ins\C$. We have the $d$-fold ramified covering
$$\rho:\Yt_S\defs{\rm Specan}_{\YS}\mopl_{k=0}^{d-1}\,\OO_{\YS}(-k)\to\YS,$$
which is totally ramified along $Z_S\sst\YS$. Here ${\rm Specan}$ means the analytic space associated with the corresponding algebraic object, and $\OO_{\YS}(-k)$ is the pullback of $\OO_Y(-k)$. The finite direct sum has an $\OO_{\YS}$-algebra structure by the morphisms $\OO_{\YS}({-}k{-}d)\tos\OO_{\YS}(-k)$ for $k\ins\N$ induced by the multiplication by $f\!_S$. It may be better to divide the $\OO_{\YS}$-algebra $\mopl_{k\in\N}\,\OO_{\YS}(-k)$ by the ideal generated locally by $w^d\mi f\!_S\one w^d$, where $w$ is a local generator of $\OO_{\YS}(-1)$ and $f\!_S$ is identified with a section of $\OO_{\YS}(d)$. So $\Yt_S$ is a divisor on the line bundle corresponding to the invertible sheaf $\OO_{\YS}(1)$, and is locally defined by
\begin{equation} \label{7.1}
w^d\mi h\mi sg'\eq0.
\end{equation}
Here $h,g'$ are respectively the restrictions of $f,g$ to $\{x_l\eq1\}$ for some $l\ins[0,n]$, and $w$ is identified with a local section of the dual line bundle, that is, a linear function on the fibers. Since $g'$ does not vanish at the singular points of $Z$, we see that $\Yt_S$ is {\it smooth\one} shrinking $S$ if necessary, since the projection $\YS\to S$ is proper.
\sk
We have the direct sum decomposition
\begin{equation} \label{7.2}
\rho_*\C_{\Yt_S}=\C_{Y_S}\oplus\mopl_{i=1}^{d-1}\,(j_S)_!L_{S,i}
\end{equation}
in a compatible way with the action of the covering transformation group $\Ga\cong\Z/d\Z$, where the $L_{S,i}$ are local systems of rank one on $U_S\defs\YS\stm Z_S$ with $j_S\col U_S\into\YS$ the inclusion such that the action on $(j_S)_!L_{S,i}$ by the generator corresponding to the inverse of the ``geometric monodromy" (see \cite{DiSa0}) is given by the multiplication by $e^{2\pi\sqrt{-1}\,i/d}$ for $i\ins[1,d{-}1]$. Let $L_{c,i}$ be the restriction of $L_{S,i}$ to $U_c\defs U_S\cap(Y{\times}\{c\})$. This can be identified up to nonzero constant multiplication with the local system $L_i$ in Section\,\,\ref{S1} for $f\!_c$. (Note that a local system of rank~1 on a Zariski-open subset of $\PP^n$ can be determined up to nonzero constant multiplication by the local monodromies, considering the tensor product of a local system with the dual of another local system, since $\PP^n$ is simply connected.)
\sk
We then get the decomposition of filtered regular holonomic $\D$-modules on a neighborhood of $Y{\times}\{0\}$\,:
\begin{equation} \label{7.3}
\rho^{\D}_*(\OO_{\Yt_S},F)=(\OO_{Y_S},F)\oplus\mopl_{i=1}^{d-1}(\M_{S,i},F)
\end{equation}
in a compatible way with the action of $\Ga$, using the {\it complete reducibility\one} of finite groups. Here $\rho^{\D}_*$ denotes the direct image as a filtered $\D$-module, and $\M_{S,i}$ is a regular holonomic $\D$-module whose dual is the regular meromorphic extension of $\OO_{U_S}{\otimes}_{\C}L_{S,d-i}$. (We can see that $L_{S,d-i}$ is the dual of $L_{S,i}$ looking at their local monodromies.) These filtered $\D$-modules are direct factors of underlying filtered $\D$-modules of mixed Hodge modules.
\sk
For $i\ins[1,d]$, $c\ins S$, put
$$\chi^p(L_{c,i})\defs\msum_{k=0}^n\,(-1)^k\dim\Gr_F^pH^k(U_c,L_{c,i}),$$
and similarly for $\chi_!^p(L_{c,i})$ with $H^k$ replaced by $H^k_c$. (Here the usual notation $\chi_c^p$ is not used to avoid a confusion related to $c\ins S$.) By duality we have
\begin{equation} \label{7.4}
\chi_!^p(L_{c,i})=(-1)^n\chi^{n-p}(L_{c,\one d-i})\q\h{for any}\,\,p\ins\Z,\,c\ins S,\,i\ins[0,d],
\end{equation}
with $L_{c,0}\eq L_{c,d}\defs\C_{U_c}$, since $L_{c,i}$ is the dual local system of $L_{c,\one d-i}$.
\sk
For $i\ins[1,d{-}1]$ and for a general $c\ins S$, we get
\begin{equation} \label{7.5}
\chi_!^p(L_{0,i})=\chi_!^p(L_{c,i})-(-1)^n\msum_{j=1}^{q''}\dim\Gr_F^p(i_j)_{\D}^*\,\varphi_s\M_{S,i},
\end{equation}
by applying the distinguished triangle $i_0\tos\psi_s\tos\varphi_s\,{\buildrel{[1]\,\,}\over\to}$ in Remark\,\,\ref{R4} to $\M_{S,i}$ and employing the isomorphism $i_0^*(j_S)_!L_{S,i}\eq(j_0)_!L_{0,i}$ together with the commutativity of the nearby and vanishing cycle functors with direct images by proper morphisms. Here $i_j\col\{p_j\}\into Y$ and $j_c\col U_c\into Y$ denote natural inclusions, and $(i_j)_{\D}^*(i_j)^{\D}_*(H,F)\eq(H,F)$ for a finite dimensional filtered $\C$-vector space. Note that the $\varphi_s\M_{S,i}$ are supported on the singular points of $Z_0\eq Z$, and is the direct sum of the direct images of finite dimensional $\C$-vector spaces by the inclusions $i_j\col\{p_j\}\into Y$, since it is a direct factor of the underlying filtered $\D$-module of a mixed Hodge module. (Note that ${\rm DR}(\M_{S,i})$ is generically a local system shifted by $n$, and this gives the sign $(-1)^n$ in \eqref{7.5}.)
\sk
In view of \eqref{7.4} and \eqref{7.5} we see that the vanishing cycle $\D$-module $\varphi_s\M_{S,i}$ gives the ``contribution of singular points". This can be calculated by using the mixed Hodge module $\varphi_s\Q_{h,\Yt_S}$ (since it commutes with the direct image by the finite morphism $\rho$) provided that we do not forget the action of the cyclic covering transformation group $\Ga\eq\Z/d\Z$, which acts on the variable $w$ in \eqref{7.2} via constant multiplications, see also \cite[(1.3.1)]{BS}. Since $g'$ does not vanish on a neighborhood of singular points of $Z$, the equation \eqref{7.1} is equivalent to 
\begin{equation} \label{7.6}
\widetilde{h}:=w'{}^d\mi h/g'=s,
\end{equation}
where $w'\defs w/g'{}^{1/d}$ locally on $Y$. The assertion for $i\ins[1,d{-}1]$ is then verified by using the Thom--Sebastiani-type theorem (see \cite{SS}, \cite{Va0}, and Remark\,\,\ref{R1.3}) together with the ``symmetry" of spectral numbers in the isolated singularity case (see Remark\,\,\ref{R1.2}). Note that $\Sp_{h''}(t)\eq\msum_{i=1}^{d-1}\,t^{i/d}$ with $h''\defs w'{}^d$, and $t^{k/d}$ corresponds in the proof of the Thom--Sebastiani-type theorem to $w'{}^{k-1}\ddd w'$, on which the covering transformation group acts as in \cite[(1.3.1)]{BS}. We have the decomposition
$$\Sp_{\widetilde{h}}(t)=\msum_{i=1}^{d-1}\,\Sp_{h/g'}(t)\,t^{i/d},$$
where the summand $\Sp_{h/g'}(t)\,t^{i/d}$ corresponds to the direct factor on which the action of the monodromy $T$ is given by multiplication by $e^{2\pi\sqrt{-1}\one i/d}$ using \cite[(1.3.1)]{BS}. The assertion for $i\ins[1,d{-}1]$ then follows using the duality \eqref{7.4}.
\sk
The argument is similar in the case $i\eq d$ using the duality \eqref{7.4} for $i\eq0,d$, and applying the above distinguished triangle of functors to (the direct image of) the constant Hodge module on $\Yt_S$ (which is non-singular shrinking $S$ if necessary), since $\chi_!^p(U_c)\eq\chi^p(Y)\mi\chi^p(Z_c)$ for any $p\ins\Z,\,c\ins S$. This finishes the proof of the first part of Theorem\,\,\ref{T1}.

\section{Proof of the second part of Theorem~\ref{T1}} \label{S8}
Set $X_c\defs f^{-1}(c)$, $X'_{c'}\defs f'{}^{-1}(c')\sst X\eq\C^{n+1}$ for $c,c'\ins\C^*$, and $\zeta_p\defs e^{2\pi\sqrt{-1}/p}$ for $p\ins\Z_{>0}$. We have the decomposition
\begin{equation} \label{8.1}
X_c=\h{$\bigsqcup$}_{p=0}^{m-1}\,X'_{\zeta_m^pc'}\q\h{for}\,\,\,c\eq c'{}^m.
\end{equation}
Note that $f\mi c\eq\mprod_{k=0}^{m-1}\,(f'\mi\zeta_m^kc')$. For $j\ins\Z$, set
$$H^j_c\defs H^j(X_c,\Q),\q H^{\prime j}_{c',p}\defs H^j(X'_{\zeta_m^pc'},\Q).$$
By \eqref{8.1} we have the direct sum decomposition
\begin{equation} \label{8.2}
H^j_c\eq\mopl_{p=0}^{m-1}\,H^{\prime j}_{c',p}\q\h{for}\,\,\,j\ins\Z.
\end{equation}
\sk
We have the {\it geometric monodromies\one} $\ga\ins{\rm Aut}(X_c)$ and $\ga'\ins{\rm Aut}(X'_{c'})$ of the homogeneous polynomials $f,f'$, which are defined respectively by
$$\ga(x_0,\dots,x_n)\eq(\zeta_d\one x_0,\dots,\zeta_d\one x_n)\q\h{and}\q\ga'(x_0,\dots,x_n)\eq(\zeta_{d'}x_0,\dots,\zeta_{d'}x_n),$$
hence
\begin{equation} \label{8.3}
\ga^m\eq\ga',
\end{equation}
since $\zeta_d^m\eq\zeta_{d'}$. Recall that $\ga^*$, $\ga'{}^*$ are the {\it inverses\one} of the monodromies, see Section\,\,\ref{S2} and \cite{DiSa0}. We have the isomorphisms
\begin{equation} \label{8.4}
\ga^*:H^{\prime j}_{c',p}\simto H^{\prime j}_{c',p-1}\q\h{for any}\,\,\,p\ins[1,m{-}1],
\end{equation}
which identify the $H^{\prime j}_{c',p}$ ($p\ins[0,m{-}1]$) with each other.
\sk
Define a matrix $A\eq(a_{i,j})_{i,j\in[1,n]}$ by
$$a_{i,j}\defs\begin{cases}\la&\h{if}\,\,\,i\eq m,\,j\eq 1,\\ 1&\h{if}\,\,\,j\eq i{+}1,\\ 0&\h{otherwise}.\end{cases}$$
We have the identity
\begin{equation} \label{8.5}
\det(xI\mi A)=x^m\mi\la,
\end{equation}
where $I$ denotes the identity matrix of size $m{\times}m$. So the assertion follows from \eqref{8.3}, \eqref{8.4} (restricting to one-dimensional $\la$-eigenspaces of $\ga'{}^*$). This finishes the proof of the second part of Theorem\,\,\ref{T1} in view of Remark\,\,\ref{R8.1} below. Note that $+(-1)^n$ comes from the difference between $H^{\ssb}$ and $\Ht^{\ssb}$. (Here we do not need the assumption that $Z'$ has only isolated singularities. The  second part of Theorem\,\,\ref{T1} would be known to specialists.)

\begin{rem} \label{R8.1}
Set $S\defs\C$ with {\it natural coordinate\one} $t$ so that $f^*t\eq f$. We have the cartesian diagram
$$\xymatrix{\,X_0\, \ar@{^(->}[r] \ar[d]^{a_{X_0}} & X \ar[d]^f & \,X_c\, \ar@{_(->}[l] \ar[d] \\ \,\{0\}\, \ar@{^(->}[r] & S & \,\{c\}\, \ar@{_(->}[l]}$$
together with the canonical isomorphisms of mixed Hodge structures
\begin{equation} \label{8.6}
{}^{\bf p}\psi_tH^jf_*(\Q_{h,X}[n{+}1])\simto H^j(a_{X_0})_*{}{}^{\bf p}\psi_f(\Q_{h,X}[n{+}1]).
\end{equation}
Here ${}^{\bf p}\psi_f\defs\psi_f[-1]$ (similarly for ${}^{\bf p}\psi_t$) so that mixed Hodge modules are preserved. (Actually these shifted functors are used in \cite{mhm}.) These isomorphisms are known to hold if $f$ is proper, see \cite[Theorem 2.14]{mhm} for the projective case. It is easy to construct the ``canonical" morphism using a compactification of $f$, and it is enough to show that this is bijective forgetting the mixed Hodge structure. Here we may replace $f\col X\tos S$ with its restriction to $X_a\cap S_b$ (over $S_b$) for $b\,{\ll}\,a^d$ using the $\C^*$-action (together with the Mittag-Leffler condition), where $X_a\defs\{|x|\slt a\}\sst X$, $S_b\defs\{|s|\slt b\}\sst S$ for $a,b\ins\R_{>0}$. The assertion then follows from the local triviality of $f$ along the $\C^*$-orbits.
\sk
We have moreover the canonical isomorphisms
\begin{equation} \label{8.7}
H^j(a_{X_0})_*{}{}^{\bf p}\psi_f(\Q_{h,X}[n{+}1])\simto H^ji_0^*{}{}^{\bf p}\psi_f(\Q_{h,X}[n{+}1]),
\end{equation}
since the restriction of the underlying $\Q$-complex of ${}^{\bf p}\psi_f(\Q_{h,X}[n{+}1])|_{X_0\setminus\{0\}}$ to each $\C^*$-orbit is locally constant. Here $i_0\col\{0\}\into X_0$ is a natural inclusion.
\sk
Note also that the (shifted) variation of mixed Hodge structure $H^jf_*(\Q_{h,X}[n{+}1])|_{S\setminus\{0\}}$ is monodromical. (Indeed, its pullback to a $d$-fold cyclic ramified covering is trivial using the $\C^*$-action. Hence it is a direct factor of the direct image of the shifted constant variation by this morphism.) So the left-hand side of \eqref{8.6} is naturally identified with its stalk at $t\eq1$ using the natural coordinate $t$.
\end{rem}

\part{Explicit computations} \label{Pa3}
\nin
In this part we explain some explicit computations of examples in the ordinary singularity case.

\section{Example I} \label{S9}
Assume the $f_k$ for $k\ins[1,r{-}2]$ are sufficiently general linear combinations of $x^2\mi z^2$ and $y^2\mi z^2$, and $f_{r-1}\eq x{-}y$, $f_r\eq x{+}y$. Then the reduced curve $\Cr$ has four singular points with reduced multiplicity $r{-}1$. Assume $C_1$ and $C_r$ have multiplicities $a$ and $b$ respectively, and the others have 1. Assume for instance $a\eq2$, $b\eq5$, $r\eq5$ so that $d\eq 14, d'\eq8$.
$$\h{$\setlength{\unitlength}{10mm}
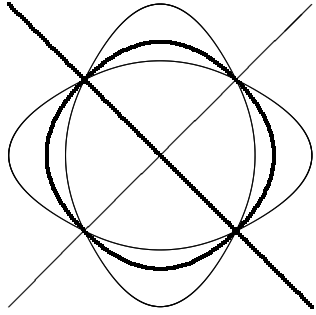
\begin{picture}(4.3,4.3)(1,1)
\put(1,1){\line(1,1){4}}
\qbezier(2,2)(3,1.5)(4,2)
\qbezier(4,2)(6,3)(4,4)
\qbezier(4,4)(3,4.5)(2,4)
\qbezier(2,4)(0,3)(2,2)
\qbezier(2,2)(3,0)(4,2)
\qbezier(4,2)(4.5,3)(4,4)
\qbezier(4,4)(3,6)(2,4)
\qbezier(2,4)(1.5,3)(2,2)
\linethickness{.5mm}
\qbezier(2,2)(3,1)(4,2)
\qbezier(4,2)(5,3)(4,4)
\qbezier(4,4)(3,5)(2,4)
\qbezier(2,4)(1,3)(2,2)
\qbezier(1,5)(3,3)(5,1)
\end{picture}$}$$
We can apply the following code of Singular \cite{Sing} which calculates the multiplicities of spectral numbers $n_{f,\frac{i}{d}}$ for $i\ins[1,3d{-}1]$ and for $a,b\ins[1,5],c\eq r\mi3\ins[0,4]$ (where one may have to change the definition of {\smaller\sf\verb@o@} in the second line if one wants to replace 5 in the second and third lines with larger numbers):
\ms
\vbox{\fontsize{8pt}{4mm}\sf\pv@ring S=0,(x),ds; int a,b,c,d,dr,dsq,i,io,j,k,l,m,n,o,p,q,r,s,w,OD; intvec GlCmp,Si,LG;@
\pv@number ga,u,v; list as,ds,z; o=20; intmat al[o][2*o]; for (a=1; a<=5; a++) {for (b=1;@
\pv@b<=5; b++) {for (c=0; c<5; c++) {GlCmp=-1,2,a,-c,2,1,-1,1,b,-1,1,1; Si=-2,c+2,a,b,-2,@
\pv@c+2,a; OD=1; LG=0; s=size(LG) div 2; w=0; for (k=1; k<=s; k++) {w=w+LG[2*k-1]*LG[2*k]*@
\pv@(LG[2*k]-1)div 2;} dsq=w; p=size(GlCmp)div 3; r=0; for(k=1; k<=p; k++) {for (i=1; i<=@
\pv@-GlCmp[3*k-2]; i++) {r++; ds[r]=GlCmp[3*k-1]; as[r]=GlCmp[3*k];}} m=100; for(k=1; k<=r;@
\pv@k++) {dsq=dsq+ds[k]*(ds[k]-1) div 2;} d=0; dr=0; for(k=1; k<=r; k++) {d=d+ds[k]*as[k];@
\pv@dr=dr+ds[k];} intmat sp[4][d]; n=size(Si); for (i=1;i<=n; i++) {z[i]=Si[i];} z[n+1]=0;@
\pv@j=0; l=1; p=1; while (p<n) {j++; s=-z[p]; al[j,1]=z[p+1];for (l=al[j,1]+1; l>1; l--)@
\pv@{al[j,l]=1;} p=p+2; for (; p<=n && z[p]>0; p++) {l++; al[j,l]=z[p];} for (i=1; i<s;@
\pv@i++) {for (l=1; l<=al[j,1]+1; l++) {al[j+i,l]=al[j,l];}} j=j+s-1;}q=j; u=1; for (i=1;@
\pv@i<=d; i++) {s=0; for (k=1; k<=r; k++) {s=s+ds[k]*(m-int(m-u*as[k]*i/d)-1);} io=i-s;@
\pv@sp[1,i]=(io-1)*(io-2)div 2; sp[2,i]=dsq+(io-1)*(dr-io-1);sp[3,i]=(dr-io-1)*(dr-io-2)@
\pv@div 2; for (j=1; j<=q; j++) {ga=0; for (l=2; l<=al[j,1]+1;l++) {v=u*al[j,l]*i/d; ga=@
\pv@ga+v-(m-int(m-v))+1;} p=m-int(m-ga); sp[1,i]=sp[1,i]-(p-1)*(p-2)div 2; sp[2,i]=sp[2,i]@
\pv@-(p-1)*(al[j,1]-p); sp[3,i]=sp[3,i]-(al[j,1]-p)*(al[j,1]-p-1)div 2;}}for (i=1; i<=d;@
\pv@i++) {for (l=1; l<=3; l++) {sp[4,i]=sp[4,i]+sp[l,i];}} p=0; for (j=1; j<=q; j++) {p=@
\pv@p+(al[j,1]-1)^2;} p=dr*(dr-3)+3-p-OD; sprintf("(a,b,c)=(
\pv@sprintf("chi(U)=
\msn
The degrees and multiplicities of global components are given by {\smaller\sf\verb@GlCmp@}, where the number of components having the same degree and multiplicity is noted together with the {\it minus sign\one} before these two numbers. (It must be written always even if it is {\smaller\sf\verb@-1@}. The minus sign is a kind of {\it separator.}) The number of local irreducible components with their multiplicities $a_{jl}$ at each singular point $p_j$ (with reduced multiplicities at least 3) is given by {\smaller\sf\verb@Si@}, where the local multiplicity can be {\it omitted\one} if it is 1 (since the number of local irreducible components is known), and the number of singular points having the {\it same data\one} is noted together with the minus sign before the data. The order of local multiplicities is irrelevant, and they may be noted weakly decreasingly. One can see the completed list of the $a_{j,l}$ by typing ``{\smaller\sf\verb@al;@}". Finally the number of ordinary points which are not contained in {\smaller\sf\verb@Si@} must be noted in {\smaller\sf\verb@OD@}. (It is not necessary to change the definition of {\smaller\sf\verb@LG@}.) The numbers of global components and singular points (with reduced multiplicities at least 3) must be at most 10, and every global reduced irreducible component must be smooth. Running this code, we get the multiplicities of spectral numbers $n_{f,\frac{i}{d}}$ for $i\ins[1,3d{-}1]$ with $(a,b,c)\eq(2,5,2)$ as follows.
$$\fontsize{10pt}{4mm}\aligned&0,0,0,1,1,1,2,1,1,2,2,2,3,9,\\&3,4,4,3,4,4,3,4,4,3,4,4,3,\h{-}4,\\&3,2,2,2,1,1,1,1,1,1,0,0,0.\endaligned$$
It computes also $\chi(U)$ and $n_{f,\frac{i}{d}}\pl n_{f,\frac{i}{d}+1}\pl n_{f,\frac{i}{d}+2}$ for $i\ins[1,d]$. These coincide with $6$ putting $n_{f,3}\eq1$. (Note that $(d'{-}3)d'{+}3\eq43$, $\nu_4\eq4$, $\nu_2\eq1$.) It seems that $n_{f,\frac{3}{d}}$ vanishes when $b\sgt a\pl c$ or $ a\sgt b\pl 2c\pl 1$.

\section{Example II} \label{S10}
Assume $f_k\eq g_kz\pl h_k$ with $g_k,h_k$ sufficiently general polynomials in $x,y$ of degree $e{-}1$ and $e$ respectively for $k\ins[1,r{-}2]$ with $e\ins\Z_{\ges2}$, and $f_{r-1}\eq x$, $f_r\eq y$. Assume the multiplicities of $C_1$ and $C_r$ are $a$ and $b$ respectively, and the others have 1. If we assume for instance $r\eq5$, $a\eq3$, $b\eq2$, $e\eq3$, then $d\eq18,d'\eq11$. The multiplicities of spectral numbers $n_{f,\frac{i}{d}}$ for $i\ins[1,3d{-}1]$ with $a,b\ins[1,5],c\eq r\mi3\ins[0,4]$ can be calculated by the code in Section\,\,\ref{S9} replacing the definitions of {\smaller\sf\verb@GlCmp@}, {\smaller\sf\verb@Si@}, {\smaller\sf\verb@OD@} and {\smaller\sf\verb@LG@} in the third and fourth lines with the following:
\ms
\vbox{\fontsize{8pt}{4mm}\sf\pv@GlCmp=-1,3,a,-c,3,1,-1,1,b,-1,1,1; Si=-1,4+2*c,a,a,b; OD=5*c*(c+1)div 2+2*(c+1); LG=@
\pv@-(c+1),2,-2,1;@}
\msn
Here {\smaller\sf\verb@LG@} is the list of the numbers $m_{j,k}$, where the multiplicity of each number is noted before the number together with the minus sign as in Section\,\,\ref{S11}.
The order of numbers is irrelevant, since we need only the sum of the $\binom{m_{j,k}}{2}$ (which can be seen up to sign by typing {\smaller\sf\verb@w;@}), see Remark\,\,\ref{R6.1}. Running the modified code, we get for $(a,b,c)\eq(3,2,2)$ the following:
$$\fontsize{10pt}{4mm}\aligned&0,0,0,2,3,4,0,3,4,4,9,11,4,5,11,13,15,24,\\&6,8,10,14,14,14,12,14,14,14,12,10,14,14,10,8,6,-4,\\&15,13,11,5,4,3,9,4,3,3,0,0,3,2,0,0,0.
\endaligned$$
Here $\chi(U)\eq21$ by Remark\,\,\ref{R4}, since $(d'{-}3)d'{+}3\eq91$, $\nu_2\eq(9{-}4){\cdot}\one 3\pl 2\one\one {\cdot}\one 3\eq21$, and $\nu_8\eq1$. It seems that $n_{f,\frac{3}{d}}$ vanishes for any $a,b\gess1$ and $c\eq r\mi3\gess0$. Note that $n_{f,\frac{i}{d}}$ for $i\ins[1,d{-}1]$ quite often vanishes when $c\eq0$ (especially if $a\eq b\eq1$).

\begin{rem} \label{R10.1}
If we assume $e\eq4$ instead of 3 in the above example, then the above code should be replaced by
\ms
\vbox{\fontsize{8pt}{4mm}\sf\pv@GlCmp=-1,4,a,-c,4,1,-1,1,b,-1,1,1; Si=-1,5+3*c,a,a,a,b; OD=7*c*(c+1)div 2+2*(c+1); LG=@
\pv@-(c+1),3,-2,1;@}
\msn
It seems that $n_{f,\frac{3}{d}}\eq0$ for any $a,b\gess1$ and $c\eq r\mi3\gess0$.
\end{rem}

\section{Example III} \label{S11}
Assume the $f_k$ are sufficiently general linear combinations of $(x^2{-}z^2)^3$ and $(y^2{-}z^2)^3$ for $k\ins[1,r{-}1]$, and $f_r\eq x\mi y$. Assume the multiplicities of $C_1$ and $C_2$ are $a$ and $b$ respectively, and 1 otherwise. Replacing the definitions of {\smaller\sf\verb@GlCmp@}, {\smaller\sf\verb@Si@}, {\smaller\sf\verb@OD@} and {\smaller\sf\verb@LG@} in the third and fourth lines of the code in Section\,\,\ref{S9} with
\ms
\vbox{\fontsize{8pt}{4mm}\sf\pv@GlCmp=-1,6,a,-1,6,b,-c,6,1,-1,1,1; Si=-2,7+3*c,a,a,a,b,b,b,-2,6+3*c,a,a,a,b,b,b; OD=0; @
\pv@LG=-4*(c+2),3,-2,1;@}
\msn
we get for $r\eq c\pl3\eq4$, $a\eq3$, $b\eq2$ the following.
$$\fontsize{9pt}{4mm}\aligned&0,0,1,1,2,2,3,3,4,4,5,5,3,3,4,4,5,5,3,3,4,4,5,5,3,3,4,4,5,5,6,6,7,7,8,8,25,\\&9,9,9,9,9,9,9,9,9,9,9,9,9,9,9,9,9,9,9,9,9,9,9,9,9,9,9,9,9,9,9,9,9,9,9,9,-9,\\&8,8,7,7,6,6,5,5,4,4,3,3,5,5,4,4,3,3,5,5,4,4,3,3,5,5,4,4,3,3,2,2,1,1,0,0.\endaligned$$
Here $d\eq37$, $d'\eq19$, $d'(d'{-}3)\pl3=307$, $\nu_{10}\eq\nu_9\eq2$, and $2(8^2\pl9^2)\eq290$. It seems that the $n_{f,\frac{i}{d}+1}$ are independent of $i\ins[1,d{-}1]$ and equal to the half of $\chi(U){+}1$ for any $a,b\gess1$, $c\eq r\mi3\gess0$. (There seem to be certain cancellations between global and local terms.)

\section{Example IV} \label{S12}
Assume $f\eq x^ay^bz^{c+1}(x\pl y)(x\pl y\pl z)$, where $a,b,c{+}1\gess 1$. (This is actually equivalent to that $f\eq x^ay^{c+1}z(x{-}z)^b(y{-}z)$ up to sign.) Replacing the definitions of {\smaller\sf\verb@GlCmp@}, {\smaller\sf\verb@Si@}, {\smaller\sf\verb@OD@} and {\smaller\sf\verb@LG@} in the third and fourth lines of the code in Section\,\,\ref{S9} with
\ms
\vbox{\fontsize{8pt}{4mm}\sf\pv@GlCmp=-1,1,a,-1,1,b,-1,1,c+1,-2,1,1; Si=-1,3,a,b,-1,3,c+1; OD=4; LG=-6,1;@} 
\msn
we get for $(a,b,c{+}1)\eq(4,2,1)$ and $(5,1,1)$ respectively the following.
$$\fontsize{10pt}{4mm}\aligned0,0,0,0,0,1,0,1,4,\,\,&\q\q\q\q0,0,0,0,1,1,1,1,4,\\
0,1,1,1,1,0,1,0,\h{-}4,&\q\q\q\q0,0,1,0,0,0,0,0,\h{-}4,\\
1,0,0,0,0,0,0,0,\q\,&\q\q\q\q1,1,0,1,0,0,0,0.\endaligned$$
These two have the same pole order spectral sequence and $^P\!n_{f,\frac{3}{d}}\eq0$, see \cite[Example 5.5]{nwh}. It does not seem easy to find a simple rule about $(a,b,c)$ for vanishing of $n_{f,\frac{3}{d}}$ although it seems to be given by $a{+}b\gess4$ when $c\eq0$.

\begin{rem} \label{R12.1}
We can also consider the case $f\eq x^ay^bz^{c+1}(xy\pl yz\pl xz)$. Here the above code is replaced by 
\ms
\vbox{\fontsize{8pt}{4mm}\sf\pv@GlCmp=-1,1,a,-1,1,b,-1,1,c+1,-1,2,1; Si=-1,3,a,c+1,-1,3,b,c+1,-1,3,a,b; OD=0; LG=-9,1;@}
\msn
For $(a,b,c{+}1)\eq(1,1,5)$ or $(1,5,1)$ or $(5,1,1)$, we get
$$\fontsize{10pt}{4mm}\aligned&0,0,0,0,1,1,0,0,3,\\&1,1,1,0,0,0,1,1,\h{-}3,\\&0,0,0,1,0,0,0,0,\endaligned$$
where the pole order spectral sequence is the same as in \cite[Example 5.5]{nwh} except that $\mu^{(r)}_9$ and $\nu^{(r)}_{18}$ for $r\eq2,3$ are changed from 4 to 3, in particular, $^P\!n_{f,\frac{3}{d}}\eq0$ as in the case of the above example.
\sk
There are very few examples such that $^P\!n_{f,\frac{3}{d}}\eq0$ and $\chi(U)\nes0$ in the ordinary singularity case, where the second condition is needed for the monodromy conjecture. A necessary condition for the first seems to be that there are $p,p'\ins C$ with multiplicities respectively $m$ and $m'$ satisfying $\tfrac{2}{m}\slt\tfrac{3}{d}$, $\tfrac{3}{d}\ins\tfrac{1}{m'}\Z$ and ${\rm mult}_{p'}\Cr\gess3$ as far as calculated (but these are never sufficient). The situation seems quite compatible with \cite[Remark 4.2c]{wh} about the failure of the converse of \cite[Theorem 2\,(a)]{bCM}. In the line arrangement case, the above condition is satisfied in the case $d\eq 3e$ for $e\ins\Z$ with $m\sgt 2e$, $m'\eq e$. Modifying a small computer program used in \cite{nwh} and calculating only the $(3{+}jd)$-degree part with $j\eq0,1,2$, we get for $d\eq 3e$ with $e\eq 4,5,6$, the vanishing of $^P\!n_{f,\frac{3}{d}}$ (with $\chi(U)\eq1$ or 2) if
$$\aligned f\eq x^{e-2}yz^{2e-1}(xy{+}yz{+}xz)\q\h{or}\q x^{e-3}yz^{2e-1}(x{-}y)(xy{+}yz{+}xz)\\ \h{or}\q x^{2e-1}y^{e-2}z(x{-}z)(y{-}z)\q\h{or}\q x^{2e-1}y^{e-3}z(x{-}z)(y^2{-}z^2),\endaligned$$
but not if $f\eq yz^5(x{-}y)(xy{+}yz{+}xz)$ or $x^4yz(x{-}y)(xy{+}yz{+}xz)$ or $x^4yz(y{-}z)(x^2{-}z^2)$. We have recently found a new example with $^P\!n_{f,\frac{3}{d}}\eq0$ for $d\eq9$ : $f\eq xz^5(x^2z{+}y^2z{+}x^3{+}xy^2)$ with $\chi(U)\eq1$. It seems quite difficult to get further examples with $^P\!n_{f,\frac{3}{d}}\eq0$ for $d\eq9$.
\end{rem}

{\makeatletter
\def\section#1#2{}%
\def\@startsection#1#2#3#4#5#6{}%
\def\@seccntformat#1{}%
\def\refname{}%
\makeatother

\par\addvspace{3ex plus 1ex minus .2ex}
\noindent{\normalfont\normalsize\bfseries References}
\addcontentsline{toc}{part}{References}  
\par\vskip 2ex plus .2ex

}
\end{document}